\documentclass[12pt]{article}
\usepackage{amssymb}
\usepackage{graphicx}
\usepackage{amsfonts}
\usepackage{latexsym}
\usepackage{amsthm}
\usepackage{amsmath}
\usepackage[all,arc]{xy} 

\usepackage{amssymb, amscd}

\usepackage{url}

\newtheorem{problem}{Problem}
\newtheorem{theo}[problem]{Theorem}

\newtheorem{defin}[problem]{Definition}
\newtheorem{prop}[problem]{Proposition}

\newtheorem{lema}[problem]{Lemma}

\begin{document}
\date{May 22, 2014}
 \title{{Symmetric polyomino tilings, \\ tribones, ideals, and Gr\" obner bases}}

\author{{Manuela Muzika Dizdarevi\' c}\\ {\small Faculty of Natural Sciences}\\[-2mm] {\small and Mathematics, Sarajevo}
\\[-2mm]{\small mmuzika@me.com}
  \and Rade  T.\ \v Zivaljevi\' c\\ {\small Mathematical Institute}\\[-2mm] {\small SASA, Belgrade}\\[-2mm]
 {\small rade$@$mi.sanu.ac.rs} }

\maketitle
\begin{abstract}
We apply the theory of Gr\" obner bases to the study of signed,
symmetric polyomino tilings of planar domains. Complementing the
results of Conway and Lagarias we show that the triangular regions
$T_N=T_{3k-1}$ and $T_N=T_{3k}$ in a hexagonal lattice admit a
{\em signed tiling} by three-in-line polyominoes (tribones) {\em
symmetric} with respect to the $120^{\circ}$ rotation of the
triangle if and only if either $N=27r-1$ or $N=27r$ for some
integer $r\geq 0$.
\end{abstract}

\renewcommand{\thefootnote}{\fnsymbol{footnote}}
\footnotetext{R.~\v Zivaljevi\'c was supported by the Grants
174017 and 174020 of the Ministry for Science and Technological
Development of Serbia.}


\section{Introduction and a summary of main results}\label{sec:main}

Our general objective is to explore {\em signed polyomino tilings}
which are symmetric with respect to a group of symmetries by the
methods of standard (Gr\" obner) bases of polynomial ideals.

\begin{figure}[hbt]
\centering
\includegraphics[scale=0.60]{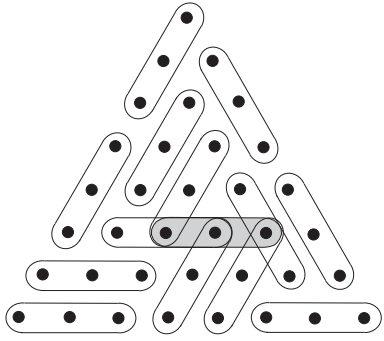}
\caption{A signed tiling of a triangular region in the dual of the
hexagonal lattice.} \label{fig:trimino-4}
\end{figure}

\medskip
The tiling depicted in Figure~\ref{fig:trimino-4}, illustrating
the case $N=8$ of Theorem~\ref{thm:Conv-Laga}, shows that a
triangular region in a hexagonal lattice may have a signed tiling
by congruent copies of the three-in-line tile (tribone). In the
same paper Conway and Lagarias showed
(\cite[Theorem~1.2.]{ConLag}) that neither this nor any other
triangular region in the hexagonal lattice can be tiled by
tribones (if `negative' tiles are not permitted).

A very nice exposition of these and related results can be found
in \cite{Thur} and \cite[Chapter 23]{FT}.

\begin{theo}\label{thm:Conv-Laga}
{\rm (Conway-Lagarias \cite[Theorem 11.4]{ConLag})} The triangular
region $T_N$ in the hexagonal lattice has a signed tiling by
congruent copies by three-in-line tiles (tribones) if and only if
$N=9r$ or $N=9r+8$ for some integer $r\geq 0$.
\end{theo}

\medskip
Our main results (Theorems~\ref{thm:main-1} and \ref{thm:main-2})
say that the triangular regions $T_N=T_{3k-1}$ and $T_N=T_{3k}$ in
a hexagonal lattice\footnotemark\ admit a tiling by tribones {\em
symmetric} with respect to the rotation of the triangle through
the angle of $120^{\circ}$ degrees if and only if either $N=27r-1$
or $N=27r$ for some integer $r\geq 0$. In particular the triangle
depicted in Figure~\ref{fig:trimino-4} does not admit such a
tiling.  \footnotetext{In the remaining case $N=3k+1$ there is a
hexagon in $T_N$ fixed by the $120^\circ$-degrees rotation.}

\medskip\noindent
The method applied for the proof of this theorem is based on the
observation that the {\em tile homology group} introduced by
Conway and Lagarias in \cite{ConLag} (see also Reid \cite[Section
2]{Reid}) is naturally a module over the group ring of the
associated group of translations. This group ring is a quotient of
a polynomial ring which allows us to reduce the tiling problem to
the `submodule membership problem' and apply the theory of
Gr\"{o}bner bases.

\medskip
Surprisingly enough there are very few applications of the
algebraic method based on Gr\"{o}bner basis to problems of tilings
and tessellations and the only reference we are aware of is the
paper by Bodini and Nouvel \cite{BN}. The fact that the `tile
homology group' in the sense of \cite{Reid} is a module over a
polynomial ring offers some obvious technical advantages. One of
our objectives is to advertise this approach in the context of
signed tilings with symmetries. These problems seem to be
particularly well adopted to the algebraic method in light of the
fruitful relationship between the theory of Gr\"{obner} bases and
the theory of invariants of group actions \cite{CLO1},
\cite{Sturm}.

\section{Generalities about lattice tilings}\label{sec:intro-A}

\begin{figure}[hbt]
\centering
\includegraphics[scale=0.60]{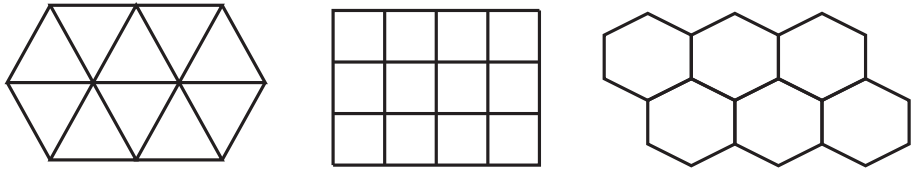}
\caption{Three regular lattice tilings of the plane.}
\label{fig:regular}
\end{figure}

There are three regular lattice tilings of $\mathbb{R}^2$,  the
triangular lattice $L_\Delta$, square lattice $L_\Box$, and the
hexagonal lattice $L_{\rm hex}$, depicted in
Figure~\ref{fig:regular}. If $L$ is one of these lattice tilings
the associated dual lattice (point set) $L^\circ$ is generated by
all barycenters of the elementary cells of $L$.

Let $A(L)$ be the free abelian group generated by all elementary
cells of the lattice $L$. A `lattice tile' $P$ (informally a
lattice figure in $L$),  defined as a finite collection
$P=\{c_1,\ldots, c_n\}$ of cells in $L$, is associated an element
$P=c_1+\ldots + c_n$ of the group $A(L)$.

\medskip
The problem if a given bounded lattice region (lattice figure) $R$
admits a {\em signed}\/ tiling with translates of prototiles
$\mathcal{R}=\{R_1,\ldots , R_k\}$ is an instance of the {\em
subgroup membership problem}. Indeed, let $B({\mathcal{R}})$ be
the subgroup of $A(L)$ generated by all translates of prototiles
$R_i$ and let $H({\mathcal{R}})=A(L)/B({\mathcal{R}})$ be the
associated `homology group'. Then (following \cite{ConLag} and
\cite{Reid}) such a tiling exists if and only if $R\in
B(\mathcal{R})$ or equivalently if the coset $R+B(\mathcal{R})$ is
the zero element in $H(\mathcal{R})$.

\medskip
Let $G = G(L)$ be the group of all affine transformations that
keep the lattice tiling $L$ invariant. Let $\Gamma=\Gamma(L)$ be
its subgroup of all translations with this property. By selecting
$0\in L^\circ$ as the zero element, $L^\circ$ is turned into a
group and there is a natural identification $\Gamma=L^\circ$.

\medskip
The group $A(L)$ is a clearly a module over the group ring
$\mathbb{Z}[\Gamma]$ (which is isomorphic to the ring
$\mathbb{Z}[\mathbb{Z}^2]$ of Laurent polynomials in two
variables. This ring can be obtained (in many ways) as a quotient
of the semigroup ring $\mathbb{Z}[\mathbb{N}^d]\cong
\mathbb{Z}[x_1,\ldots, x_d]$ (for some $d$).

\medskip
This observation allows us to see the groups $A(L),
B(\mathcal{R})$ and $H(\mathcal{R})$ as modules over the
polynomial ring $\mathbb{Z}[x_1,\ldots, x_d]$ and to reduce the
tiling question to the {\em submodule membership problem}
\cite[Chapter 5]{CLO2}. In turn, in the spirit of \cite{BN}, one
can use the ideas and methods of Gr\" obner basis theory.

\medskip
Here we put some emphasis on the use of the `submodule membership
problem' as a natural extension of the `ideal membership problem',
originally proposed and used by Bodini and Nouvel in \cite{BN}.
This appears to be a more natural and conceptual approach to the
general tiling problems since the module $A(L)$ is no longer
required to be monogenic (cyclic over $\mathbb{Z}[P]$) which
allows us a greater freedom in choosing the semigroup ring
$\mathbb{Z}[P]$. This property will be indispensable in the study
of tilings {\em symmetric} with respect to a group of symmetries
which is the main goal of this paper.

\subsection{An example}\label{sec:intro-B}

\begin{figure}[hbt]
\centering
\includegraphics[scale=0.60]{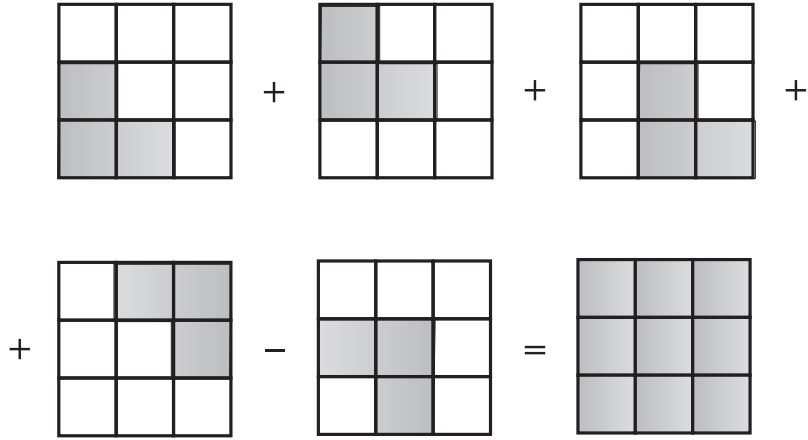}
\caption{A signed tiling of a $3\times 3$ square symmetric w.r.t.\
the main diagonal.} \label{fig:mreza-1}
\end{figure}

The reader may find the following example, depicted in
Figure~\ref{fig:mreza-1}, as a good illustration of the main
problem studied in our paper.

\medskip
The $(3\times 3)$ checkerboard $C_{3\times 3}$
(Figure~\ref{fig:mreza-1}) is supposed to be paved by translates
of two types of {\em prototiles}. Each of the cells (elementary
squares) is labelled (coordinatized) by a pair $(i,j)\in
\mathbb{N}^2$ of integers and each tile (polyomino) is formally a
union of a finite number of elementary cells. In the example
depicted in Figure~\ref{fig:mreza-1} there are two types of
prototiles, $T_1 = \{(0,0),(1,0),(0,1)\}$ and $T_2=\{(1,1),
1,0),(0,1)\}$.

The tiling depicted in Figure~\ref{fig:mreza-1} satisfies the
following conditions:

\begin{enumerate}
 \item[(1)] The $(3\times 3)$ chessboard $C_{3\times 3}$ is covered
by $3$ translated copies of prototile $T_1$ and $2$ translated
copies of prototile $T_2$. In general the translated prototiles
are not required to be subsets of $C_{3\times 3}$.
 \item[(2)] Each (translated) prototile is associated a weight (sign)
and the total weight of each of the cells $(i,j)\in C_{3\times 3}$
is equal to $1$. (The condition that the total weight of cells
outside $C_{3\times 3}$ is zero is added if prototiles are not
necessarily subsets of $C_{3\times 3}$.)
 \item[(3)] The tiling is symmetric with respect to the main diagonal
of $C_{3\times 3}$ in the sense that if a translated prototile $T
+ v$ (where $T\in\{T_1,T_2\}$ and $v\in \mathbb{N}^2$) appears in
the decomposition (tiling) with a weight $w\in \mathbb{Z}$ then
the diagonally symmetric prototile $T'+v'$ is also present with
the same weight.
\end{enumerate}
A basic observation is that each polyomino $P\subset \mathbb{N}^2$
can be associated a polynomial $f_P = \sum\{x^iy^j\mid (i,j)\in
P\}$, for example $f_{T_1}=1+x+y$ and $f_{T_2}=x+y+xy$. The
decomposition depicted in Figure~\ref{fig:mreza-1} naturally
corresponds to the following decomposition of polynomials in the
ring $\mathbb{Z}[x,y]$ (or in the ring
$\mathbb{Z}[\sigma_1,\sigma_2]$ of symmetric polynomials):
\begin{equation}\label{eqn:intro-example}
f_{T_1}+[xf_{T_1} + yf_{T_1}] - f_{T_2} + xyf_{T_2} =
(1+x+x^2)(1+y+y^2) = f_{C_{3\times 3}}.
\end{equation}
Our immediate objective is to use the theory of Gr\" obner bases
to generate such identities. More generally we want to develop and
study procedures and algorithms  for the systematic analysis
(existence and other properties) of decompositions similar to
(\ref{eqn:intro-example}).

\subsection{Basic facts about polyomino tilings}
\label{sec:vocabulary}

Informally a polyomino pattern $P$ (or polyomino for short) is a
(not necessarily connected) finite region consisting of cells in
one of the three regular lattice tilings of the plane
(Figure~\ref{fig:regular}). It is sometimes more convenient to
describe a polyomino as a collection $P=\{c_1,\ldots, c_k\}$ of
elementary cells in the associated lattice $L$. Some authors use
the generic name {\em polyforms}\/ for all three types of
polyominoes so the $L_\square$-polyforms are polyominoes in the
usual sense \cite{Gol}, $L_{hex}$-polyforms are referred to as
{\em polyhexes} etc.

\medskip
We frequently use a slightly more general (algebraic) definition
of a polyomino as a {\em multiset}, subset of $L$, with multiple
and possibly with negative elements. We will tacitly make a
distinction between the geometric and algebraic definition by
reserving the term `weighted polyomino' for the algebraic version.
However most of the time the term `polyomino' is used
interchangeably for both kinds of polyomino patterns.

\begin{defin}\label{def:polyomino-as-multiset}
A (weighted) polyomino $P$ is a finite weighted subset of $L$ (a
multiset) which contains each elementary cell $c\in L$ with some
(positive or negative) multiplicity $w_c\in \mathbb{Z}$. In other
words $P = \sum w_c c$ is an element of the free abelian group
$A(L)$ generated by all cells of the lattice tiling $L$.
\end{defin}

\subsection{Geometric-algebraic dictionary of polyomino tilings}

We have already seen in Section~\ref{sec:intro-B} an example of
the correspondence between a geometric image
(Figure~\ref{fig:mreza-1}) and an algebraic expression (equation
(\ref{eqn:intro-example})), based on the correspondence
$(i,j)\leftrightarrow x^iy^j$ between the cell labelled by
$(i,j)\in \mathbb{N}^2$ and the associated monomial $x^iy^j$.

More generally let $S\cong \mathbb{N}^d$ be a semigroup which acts
on the lattice $L$ by translations, which means that there exists
a homomorphism $\rho : S\rightarrow \Gamma$ from $S$ to the group
$\Gamma = \Gamma(L)$ of all translations that keep the lattice
tiling $L$ invariant.

The group $A(L)$ is naturally a module over the semigroup ring
$\mathbb{Z}[S]\cong \mathbb{Z}[x_1,\ldots, x_d]$. For example if
$S=\Gamma$ then $\mathbb{Z}[S]\cong \mathbb{Z}[\Gamma]$ and $A(L)$
is a $\mathbb{Z}[\Gamma]$-module where $\mathbb{Z}[\Gamma]\cong
\mathbb{Z}[x, x^{-1}; y, y^{-1}]$ is the ring of Laurent
polynomials.

Let $\mathcal{R}=\{P_1,\ldots, P_k\}$ be a collection of basic
tiles (prototiles). Define $B(\mathcal{R})$ as the subgroup of
$A(L)$ generated by all translates of the prototiles $P_i$, or
equivalently as a $\mathbb{Z}[S]$-submodule of $A(L)$ generated by
$\mathcal{R}$.

\medskip
The following tautological proposition links the idea of the tile
homology group of Conway and Lagarias \cite{ConLag} and Reid
\cite{Reid} with the `submodule membership problem' typical for
applications of Gr\"{o}bner bases (as proposed by Bodini and
Nouvel \cite{BN}).

\begin{prop}\label{prop:fund-prop}
A polyomino $P$ has a signed tiling by translates of prototiles
$\mathcal{R}=\{P_1,\ldots, P_k\}$ if and only if $P\in
B(\mathcal{R})$ where $B(\mathcal{R})$ is the
$\mathbb{Z}[S]$-submodule of $A(L)$ generated by $\mathcal{R}$.
The associated class $[P]$ in the {\em tile homology module}
\begin{equation}\label{eqn:tile-module}
H(\mathcal{R}):=A(L)/B(\mathcal{R})
\end{equation}
is a `quantitative measure' of how far is $P$ from admitting a
tiling by $\mathcal{R}$.
\end{prop}

A modified version of Proposition~\ref{prop:fund-prop} applies to
proper subsets of the lattice tiling $L$. The following
proposition serves as an illustration of the simplest case where
we restrict our attention to the first quadrant of the
$L_\square$-lattice tiling. In this case $S = \mathbb{N}^2$ and
$\mathbb{Z}[S]=\mathbb{Z}[x,y]$. As in Section~\ref{sec:intro-B},
each polyomino $P\subset \mathbb{N}^2$ is associated a polynomial
$f_P = \sum\{x^iy^j\mid (i,j)\in P\}$.

\begin{prop}\label{prop:signed} A polyomino pattern $P\subset \mathbb{N}^2$ admits a
signed tiling by the first quadrant translates of polyomino
patterns $P_1,\ldots, P_k$ if and only if, $$ f_P =
h_1f_{P_1}+\ldots + h_kf_{P_k} $$ for some polynomials
$h_1,\ldots, h_k$ with arbitrary integer coefficients or
equivalently if,
$$
f_P\in \langle f_{P_1},\ldots, f_{P_k} \rangle.
$$
\end{prop}

\subsection{Equivariant polyomino tilings}

The group $G = G(L)$ was introduced in Section~\ref{sec:intro-A}
as the group of all affine transformations that keep the lattice
tiling $L$ invariant. The abelian group $A(L)$ is a module over
the group ring $\mathbb{Z}[G]$. Since $\Gamma\subset G$ is a
normal subgroup, we observe that $G$ acts on $A(L)$ preserving its
$\mathbb{Z}[\Gamma]$-module structure as well, provided $G$ acts
on the `scalars' from $\mathbb{Z}[\Gamma]$ by conjugation.

Let $Q\subset G$ be a (finite) subgroup of $G$. Assume that the
set $\mathcal{R}$ of prototiles is invariant with respect to the
group $Q$. Then $Q$ acts on the submodule $B(\mathcal{R})$ and the
tile homology module (\ref{eqn:tile-module}). Again, one shouldn't
forget that the action of $Q$ on scalars from $\mathbb{Z}[\Gamma]$
may be non-trivial. Define $B(\mathcal{R})^Q = {\rm
Hom}_Q(\mathbb{Z}, B(\mathcal{R}))$ as the subgroup (submodule) of
$B(\mathcal{R})$ of elements which are invariant under the action
of $Q$.

\medskip
If we restrict our attention to the subring
$\mathbb{Z}[\Gamma]^Q\subset \mathbb{Z}[\Gamma]$ of $Q$-invariant
elements then the action of $Q$ on scalars from
$\mathbb{Z}[\Gamma]^Q$ is trivial and the
$\mathbb{Z}[\Gamma]^Q$-module $A(L)$ is a $Q$-module in the usual
sense.

An element of the group $A(L)^Q$ is referred to as an {\em
equivariant signed polyomino}. The fundamental problem is to
decide when a given polyomino $P\in A(L)^Q$ admits a $Q$-symmetric
signed tiling by translates of a $Q$-invariant family of
prototiles $\mathcal{R}$. The following criterion is an
equivariant analogue of Proposition~\ref{prop:fund-prop}.

\begin{prop}\label{prop:equi-fund-prop} Let $\mathcal{R}$ be a
$Q$-invariant (finite) set of prototiles.  A $Q$-invariant
polyomino $P\in A(L)^Q$ has an equivariant, signed tiling by
translates of prototiles $\mathcal{R}$ if and only if $P\in
B(\mathcal{R})^Q$ where $B(\mathcal{R})$ is the
$\mathbb{Z}[S]^Q$-submodule of $A(L)$ generated by $\mathcal{R}$.
\end{prop}
The setting of Proposition~\ref{prop:equi-fund-prop} is exactly
the same as before (Proposition~\ref{prop:fund-prop}), however the
emphasis is now on the $\mathbb{Z}[\Gamma]^Q$-module structure on
$A(L)^Q$ and $B(\mathcal{R})^Q$. In order to apply this criterion
one is supposed to determine the ring of invariants
$\mathbb{Z}[S]^Q$ and the structure of the module
$B(\mathcal{R})^Q$. Both goals can be achieved with the aid of the
theory of Gr\"{obner bases}, see \cite[Section~7]{CLO1} for
necessary tools.

\section{Hexagonal polyomino with symmetries}
\label{sec:hexagonal}

\subsection{Lattices and semigroup rings}\label{sec:hex-intro}

\begin{figure}[hbt]
\centering
\includegraphics[scale=0.50]{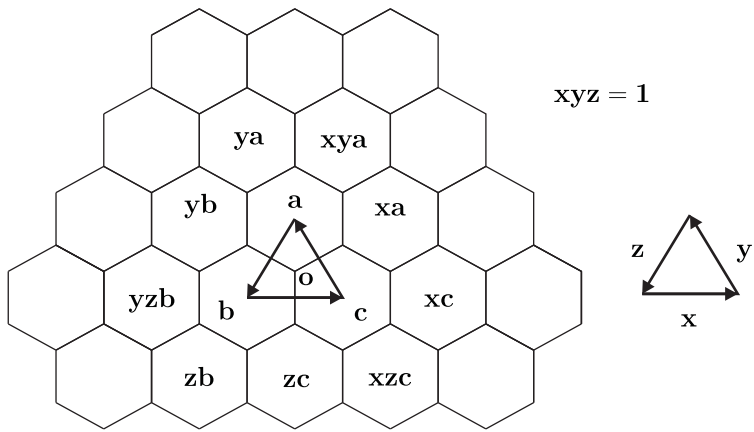}
\caption{The hexagonal tiling group $A_{hex}$ as a module over
$P=\mathbb{Z}[x,y,z]/\langle xyz-1 \rangle$.}
\label{fig:hex-mreza-4}
\end{figure}

Let $G_{hex}$ be the group of symmetries of the hexagonal tiling
$L_{hex}$ of the plane depicted in Figure~\ref{fig:hex-mreza-4}.
Our objective is to study $L_{hex}$-tiling problems which are
symmetric with respect to some (finite) subgroup of $G_{hex}$. Our
initial focus is on subgroups which act without fixed points
(invariant hexagons) so let $S_6$ be the group of all elements in
$G_{hex}$ which keep the vertex $\mathbf{o}$ fixed, and let
$\mathbb{Z}_3$ be its subgroup generated by the
$120^\circ$-rotation.


\medskip

The group $G_{hex}$ has a free abelian subgroup $D =
\Gamma(L_{hex})\cong \mathbb{Z}^2$ of rank $2$ which is generated
by three translations (vectors) $t_x, t_y, t_y$ satisfying the
condition $t_x+t_y+t_z=0$. The associated group ring $P =
\mathbb{Z}[D]$ is isomorphic to the ring $\mathbb{Z}[x,y; x^{-1},
y^{-1}]$ of Laurent polynomials in two variables. For our purposes
a more convenient representation is $P = \mathbb{Z}[x,y,z]/\langle
xyz-1 \rangle$ (Figure~\ref{fig:hex-mreza-4}) where variables
$x,y,z$ correspond to vectors $t_x = \overrightarrow{bc},
t_y=\overrightarrow{ca}, t_z = \overrightarrow{ab}$.

\begin{figure}[hbt]
\centering
\includegraphics[scale=0.90]{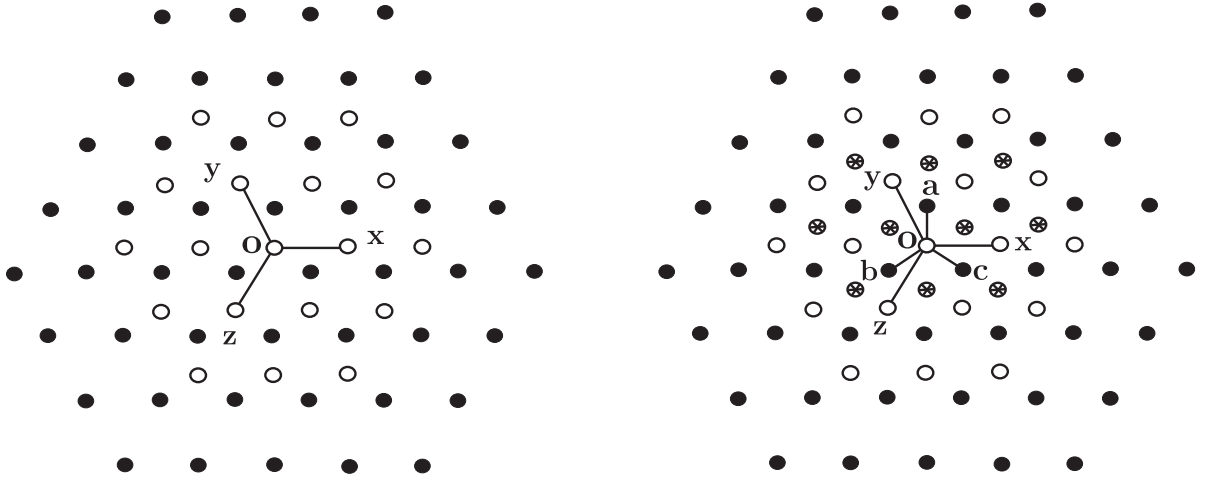}
\caption{The lattice of black dots as a module over the lattice of
white dots.} \label{fig:trougaona-mreza-2}
\end{figure}
Let $A_{hex}$ be the (infinite dimensional) free abelian group
generated by all elementary hexagonal cells of the lattice
$L_{hex}$. The group $A_{hex}$ is a finitely generated module over
the ring $P$, indeed it is generated by the three neighboring
cells $a,b,c$ with common vertex $O$, depicted in
Figure~\ref{fig:hex-mreza-4}.

For added clarity from here on the lattice $L_{hex}$ is
represented by its dual lattice $L^\circ_{hex}$ of barycenters of
all hexagons (the black dots in
Figure~\ref{fig:trougaona-mreza-2}). Consequently the `lattice'
$L^\circ_{hex}$ is a geometric object (a periodic set of points).

The lattice (discrete subgroup of $\mathbb{Z}^2$) $D$, generated
by vectors $t_x, t_y, t_z$ is in this figure represented by white
dots. The fact that $A_{hex}$ is a module over $P$ is simply a
reformulation of the fact that the lattice of white dots acts on
the set of black dots.

\medskip
The lattice $D$ is sometimes (Section~\ref{sec:abcxyz}) referred
to as the $xyz$-lattice. The lattice (discrete group) $E$
generated by black dots is referred to as the $abc$-lattice since
it is generated by vectors $t_a, t_b, t_c$, where $t_a+t_b+t_c=0$
(Figure~\ref{fig:trougaona-mreza-2}). The group ring of $E$ is $Q
= \mathbb{Z}[a,b,c]/\langle abc-1 \rangle$.

Note that $E$ has three types of points
(Figure~\ref{fig:trougaona-mreza-2} on the right) which reflects
the fact that the `white dot lattice' $D$ is a sublattice of $E$
of index $3$.

\subsection{The actions of $S_3$ and $\mathbb{Z}_3$ on $\mathbb{Z}[x,y,z]$}
\label{sec:actions-intro}

Here we collect some basic facts about the symmetric group $S_3$
and the cyclic group $\mathbb{Z}_3$ actions on $\mathbb{Z}[x,y,z]$
induced by the permutations of variables $x, y, z$. As usual for a
given $G$-module $M$ the associated submodule of $G$-invariant
elements is $M^G$. Elementary symmetric polynomials are
$\sigma_1=x+y+z, \sigma_2=xy+yz+zx, \sigma_3=xyz$.

\medskip\noindent
The $S_3$-invariant polynomials in $\mathbb{Z}[x,y,z]$ which form
a $\mathbb{Z}$-basis are $\sigma_3^p = x^py^pz^p$ (where $p\geq
0$), $\Delta(x^py^pz^q) = x^py^pz^q+y^pz^px^q+z^px^py^q$ (for
$p\neq q$), and for $p\neq q\neq r\neq p$,
\begin{equation}\label{eqn:S3-basic}
H(x^py^qz^r) = x^py^qz^r+y^pz^qx^r+z^px^qy^r+
y^px^qz^r+x^pz^qy^r+z^py^qx^r.
\end{equation}

\noindent Basic $\mathbb{Z}_3$-invariant polynomials in
$\mathbb{Z}[x,y,z]$ are,
\begin{equation}\label{eqn:Z3-basic}
 x^py^pz^p \quad (\mbox{\rm where}\, p\geq 0) \, \mbox{ {\rm and} } \,
 \Delta(x^py^qz^r) = x^py^qz^r+y^pz^qx^r+z^px^qy^r,
\end{equation}
where $(p,q,r)\neq (p,p,p)$. There is an involution $I$ on the set
$\mathbb{Z}[x,y,z]^{\mathbb{Z}_3}$ of $\mathbb{Z}_3$-invariant
polynomials defined by $I(p(x,y,z))=p(y,x,z)$. The map $\alpha :
\mathbb{Z}[x,y,z]^{S_3} \rightarrow
\mathbb{Z}[x,y,z]^{\mathbb{Z}_3}$ is a monomorphism and the image
$\mbox{\rm Im}(\alpha)$ is the fixed point set of the involution
$I$. More explicitly, $\alpha(x^py^pz^p)=x^py^pz^p$,

\begin{equation}\label{eqn:explicit}
 \alpha(\Delta(x^py^pz^q))=
\Delta(x^py^pz^q)\, \mbox{ {\rm and} }\, \alpha(H(x^py^qz^r))=
\Delta(x^py^qz^r) + I(\Delta(x^py^qz^r)).
\end{equation}

From here we deduce the following proposition.

\begin{prop}\label{prop:comm-diag-1}
There is a commutative diagram
\begin{equation}\label{eqn:comm-diag-1}
\xymatrix{ 0 \ar[r] & \langle xyz-1 \rangle^{S_3} \ar[d]^{\alpha'}
\ar[r]  &
\mathbb{Z}[x,y,z]^{S_3}\ar[d]^{\alpha}\ar[r] & \mathbb{Z}[\sigma_1,\sigma_2]\ar[r]\ar[d]^{\alpha''} & 0 \\
0 \ar[r] & \langle xyz-1 \rangle^{\mathbb{Z}_3} \ar[r] &
\mathbb{Z}[x,y,z]^{\mathbb{Z}_3} \ar[r] &
(\mathbb{Z}[x,y,z]/\langle xyz-1 \rangle)^{\mathbb{Z}_3}
  \ar[r] & 0
  }
\end{equation}
where $\langle xyz-1 \rangle\subset \mathbb{Z}[x,y,z]$ is the
principal ideal generated by $xyz-1$, with the split horizontal
exact sequences and injective vertical homomor\-phi\-sms $\alpha,
\alpha'$ and $\alpha''$.
\end{prop}

\medskip\noindent
{\bf Proof:} Since $\mathbb{Z}[x,y,z]^{S_3}=
\mathbb{Z}[\sigma_1,\sigma_2,\sigma_3]$ and $\langle xyz -
1\rangle^{S_3}= \mathbb{Z}[\sigma_1,\sigma_2,\sigma_3](\sigma_3 -
1)$, the exactness of the first row in (\ref{eqn:comm-diag-1}) is
an immediate consequence.

More explicitly the description of $S_3$-invariant and
$\mathbb{Z}_3$-invariant polynomials in $\mathbb{Z}[x,y,z]$
((\ref{eqn:S3-basic}) and (\ref{eqn:Z3-basic})) allows to describe
in a similar fashion invariant polynomials in the ideal
(submodule) $\langle xyz - 1\rangle = \langle \sigma_3-1\rangle$.
For example the basic $S_3$-invariant polynomials in $\langle
\sigma_3-1\rangle^{S_3}$ are,
\begin{equation}\label{eqn:basic-S3-inv-ideal}
\sigma_3^p(\sigma_3 - 1),\quad \Delta(x^py^pz^q)(\sigma_3 - 1),
\quad H(x^py^qz^r)(\sigma_3 - 1).
\end{equation}
\begin{figure}[htb]
\centering
\includegraphics[scale=0.60]{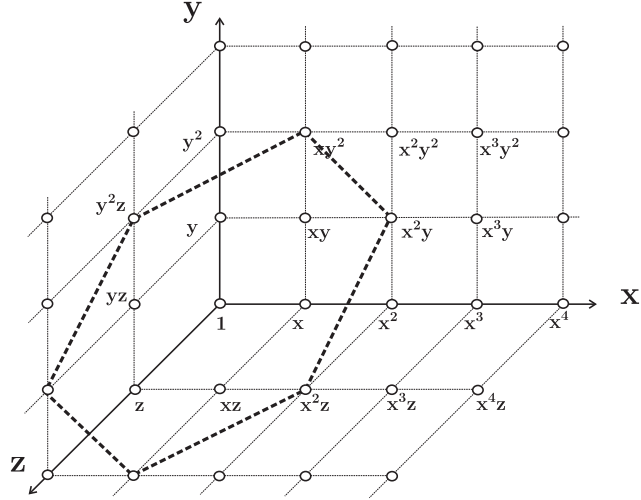}
\caption{$3d$-representation of invariant polynomials
(Proposition~\ref{prop:comm-diag-1})}
\label{fig:1-kvadrant-3d-mreza}
\end{figure}
Similarly, the basic $\mathbb{Z}_3$-invariant polynomials in
$\langle \sigma_3-1\rangle^{\mathbb{Z}_3}$ are
\begin{equation}\label{eqn:basic-Z3-inv-ideal}
\sigma_3^p(\sigma_3 - 1)\quad \mbox{\rm and} \quad
\Delta(x^py^pz^q)(\sigma_3 - 1).
\end{equation}
There is an exact sequence of $\mathbb{Z}[\mathbb{Z}_3]$-modules,
\begin{equation}
0\rightarrow \langle xyz-1 \rangle \longrightarrow
\mathbb{Z}[x,y,z] \longrightarrow \mathbb{Z}[x,y,z]/\langle xyz-1
\rangle \rightarrow 0
\end{equation}
The ideal  $\langle xyz-1 \rangle$ is as a $\mathbb{Z}$-submodule
of $\mathbb{Z}[x,y,z]$ freely generated by binomials
$x^py^qz^r(xyz-1)= x^{p+1}y^{q+1}z^{r+1}-x^py^qz^r$. This binomial
is $\mathbb{Z}_3$-invariant if and only if $p=q=r$. From here we
easily deduce the structure of $\langle xyz-1 \rangle$ as a
$\mathbb{Z}[\mathbb{Z}_3]$-module, in particular we observe that
there is a decomposition  $\langle xyz-1 \rangle\cong T\oplus F$
of $\mathbb{Z}[\mathbb{Z}_3]$-modules where $T$ is a trivial and
$F$ a free $\mathbb{Z}[\mathbb{Z}_3]$-module.

It follows that $H^1(\mathbb{Z}_3; \langle xyz-1 \rangle)\cong 0$
and from the long exact sequence of cohomology we obtain the
exactness of the second row of (\ref{eqn:comm-diag-1}),
\begin{equation}\label{eqn:z3-invariants-comm}
0\rightarrow \langle \sigma_3-1 \rangle^{\mathbb{Z}_3}
\longrightarrow \mathbb{Z}[x,y,z]^{\mathbb{Z}_3} \longrightarrow
(\mathbb{Z}[x,y,z]/\langle xyz-1 \rangle)^{\mathbb{Z}_3}
\rightarrow 0 .
\end{equation}
In particular,
\begin{equation}\label{eqn:quotient}
(\mathbb{Z}[x,y,z]/\langle xyz-1 \rangle)^{\mathbb{Z}_3} \cong
\mathbb{Z}[x,y,z]^{\mathbb{Z}_3}/\langle xyz-1
\rangle^{\mathbb{Z}_3}.
\end{equation}
The injectivity of $\alpha$ and $\alpha'$ follows from
(\ref{eqn:explicit}). In order to establish the injectivity of
$\alpha''$ we observe that (in light of (\ref{eqn:Z3-basic}) and
(\ref{eqn:basic-Z3-inv-ideal}))
$\mathbb{Z}[x,y,z]^{\mathbb{Z}_3}/\langle xyz-1
\rangle^{\mathbb{Z}_3}$ is isomorphic to the submodule of
$\mathbb{Z}[x,y,z]^{\mathbb{Z}_3}$ generated by $1 = x^0y^0z^0,
\Delta(x^p) = x^p+y^p+z^p$ and
$\Delta(x^py^q)=x^py^q+y^pz^q+z^px^q$ (where $(p,q)\neq (0,0)$).
Similarly, in light of (\ref{eqn:S3-basic}) and
(\ref{eqn:basic-S3-inv-ideal}) we observe that
$\mathbb{Z}[\sigma_1,\sigma_2]\cong
\mathbb{Z}[x,y,z]^{S_3}/\langle xyz-1 \rangle^{S_3}$ is generated
by $1 = x^0y^0z^0, \Delta(x^p), \Delta(x^py^p)$ and $H(x^py^q)$.
For example one of these $H$-polynomials (or hexagons) is depicted
in Figure~\ref{fig:1-kvadrant-3d-mreza}. From here we observe that
$\alpha''$ is injective since it satisfies analogues of formulas
(\ref{eqn:explicit}). \hfill$\square$

\subsection{The ring $P^{\mathbb{Z}_3}$ of $\mathbb{Z}_3$-invariant polynomials}
\label{sec:ring-Z3-inv}

We begin the analysis of $\mathbb{Z}_3$-invariant hexagonal
tilings by describing the structure of the ring $P^{\mathbb{Z}_3}
= (\mathbb{Z}[x,y,z]/\langle xyz-1\rangle)^{\mathbb{Z}_3}$ of
$\mathbb{Z}_3$-invariant polynomials.

\begin{figure}[htb]
\centering
\includegraphics[scale=0.80]{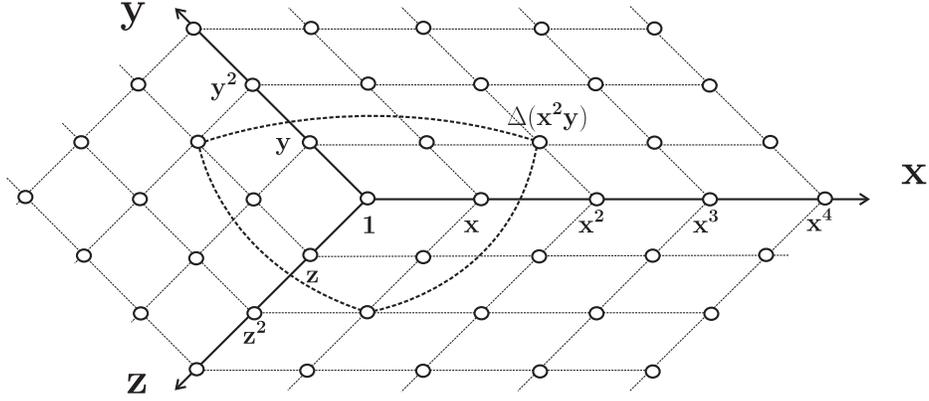}
\caption{Polynomial $\Delta(x^2y)= x^2y+y^2z+z^2y$ as one of basic
$\Delta$-polynomials which generate (over $\mathbb{Z}$) all
$\mathbb{Z}_3$-invariant polynomials in $P$.}
\label{fig:romb-mreza-2}
\end{figure}

The ring $P$ is as a $\mathbb{Z}$-module freely generated by
monomials $x^py^qz^r$ which are not divisible by $xyz$, that is
monomials $x^py^qz^r$ which satisfy at least one of conditions
$p=0, q=0, r=0$. Each of these monomials is associated a white dot
in the $XYZ$-coordinate system (Figure~\ref{fig:romb-mreza-2}).

The ring $P^{\mathbb{Z}_3}$ of $\mathbb{Z}_3$-invariant
polynomials is as a free $\mathbb{Z}$-module generated by
polynomials (`triangles') $\Delta(x^py^q) = x^py^q + y^pz^q +
z^px^q$, where either $p>0$ or $q>0$, and the constant monomial $1
= x^0y^0z^0$. Note that $\Delta(x^py^q)$ is the
$\mathbb{Z}_3$-symmetrized version of $x^py^q$; we also write
$\Delta(1)=3=1+1+1$. Moreover, $x^py^q$ is the leading monomial of
$\Delta(x^py^q)$ in a monomial order such that $x>y>z$.
Consequently, $\Delta(y^k)=\Delta(z^k)$ is almost always recorded
as $\Delta(x^k)$.

\begin{lema}\label{lema:identities}  The following identities hold in the ring $P
= \mathbb{Z}[x,y,z]/\langle xyz-1 \rangle$.
\begin{equation}\label{eqn:2-2}
\Delta(xy)\Delta(x^{p-1}y^{q-1}) = \Delta(x^py^q) +
\Delta(x^{p-1}y^{q-2})+ \Delta(x^{p-2}y^{q-1})\quad (p\geq 2,
q\geq 2)
\end{equation}
\begin{equation}\label{eqn:2-1}
\Delta(x) \Delta(x^{p-1}y^q) = \Delta(x^py^q) +
\Delta(x^{p-1}y^{q+1})+ \Delta(x^{p-2}y^{q-1})\quad \mbox{ {\rm
for} } p\geq 2 \mbox{ {\rm and} } q\geq 1
\end{equation}
\begin{equation}\label{eqn:1-2}
\Delta(x) \Delta(x^{p}y^{q-1}) = \Delta(x^py^q) +
\Delta(x^{p+1}y^{q-1})+ \Delta(x^{p-1}y^{q-2})\quad \mbox{ {\rm
for} } p\geq 1 \mbox{ {\rm and} } q\geq 2
\end{equation}
\begin{equation}\label{eqn:p>2}
\Delta(x) \Delta(x^{p-1}) = \Delta(x^p) + \Delta(x^{p-1}y)+
\Delta(xy^{p-1}) \quad \mbox{ {\rm for} } p\geq 2
\end{equation}
\begin{equation}\label{eqn:p>2bis}
\Delta(xy) \Delta(x^{p-1}) = \Delta(x^py) + \Delta(x^{p-2})+
\Delta(xy^{p}) \quad \mbox{ {\rm for} } p\geq 2.
\end{equation}
\end{lema}

By the injectivity of the map $\alpha''$
(Proposition~\ref{prop:comm-diag-1}) the ring
$\mathbb{Z}[\sigma_1, \sigma_2]$ can be seen as a subring of the
ring $P=\mathbb{Z}[x,y,z]/\langle xyz-1 \rangle$, so both $P$ and
$P^{\mathbb{Z}_3}$ are modules over $\mathbb{Z}[\sigma_1,
\sigma_2]$.

\begin{theo}\label{thm:Z3-sim-polynomials }
The ring $P^{\mathbb{Z}_3} = (\mathbb{Z}[x,y,z]/\langle xyz-1
\rangle)^{\mathbb{Z}_3}$ of\, $\mathbb{Z}_3$-invariant polynomials
is isomorphic {\rm (}as a module over
$P_\sigma=\mathbb{Z}[\sigma_1, \sigma_2]${\rm )} to the free
module $P_\sigma\cdot 1 \oplus P_\sigma\cdot\theta$ of rank two
where $\theta = \Delta(x^2y)$ {\rm
(}Figure~\ref{fig:romb-mreza-2}\,{\rm )}. Moreover,
\begin{equation}\label{eqn:theta-module} \Theta = \Theta(\sigma_1, \sigma_2, \theta)
:= \theta^2 -(\sigma_1\sigma_2 -3)\theta
+(\sigma_1^3+\sigma_2^3-6\sigma_1\sigma_2+9) = 0
\end{equation}
so there is an isomorphism of rings,
\begin{equation}\label{eqn:theta-ring}
(\mathbb{Z}[x,y,z]/\langle xyz-1 \rangle)^{\mathbb{Z}_3}\cong
\mathbb{Z}[\sigma_1,\sigma_2,\theta]/\langle \Theta \rangle
\end{equation}
where $\langle \Theta \rangle$ is the principal ideal in
$\mathbb{Z}[\sigma_1,\sigma_2,\theta]$ generated by $\Theta$.
\end{theo}

\medskip\noindent
{\bf Proof:} By definition $\Delta(x)=\sigma_1,
\Delta(xy)=\sigma_2, \theta=\Delta(x^2y)$ and
$\theta'=\Delta(xy^2)$ . Let us show that the ring
$P^{\mathbb{Z}_3}$ of $\mathbb{Z}_3$-invariant polynomials in $P$
is generated as a $\mathbb{Z}[\sigma_1, \sigma_2]$-module by $1$
and $\theta$.
\begin{itemize}
\item[$\bullet$] If $p\geq 2$ and $q\geq 2$ then by the identity
(\ref{eqn:2-2}) (Lemma~\ref{lema:identities}) the polynomial
$\Delta(x^py^q)$ can be expressed in terms of lexicographically
smaller $\Delta$-polynomials, multiplied by elements of
$\mathbb{Z}[\sigma_1, \sigma_2]$.

\item[$\bullet$] If $p\geq 3$ and $q=1$ then by
Lemma~\ref{lema:identities} (equation (\ref{eqn:2-1})) the
polynomial $\Delta(x^py)$ can be also reduced to lexicographically
smaller $\Delta$-polynomials.

\item[$\bullet$] If $p\geq 2$ then $\Delta(x^p)$ is by
Lemma~\ref{lema:identities} (equation (\ref{eqn:p>2})) reducible
to lexicographically smaller $\Delta$-polynomials.
\end{itemize}
It follows that all $\Delta$-polynomials can be expressed in terms
of $\theta, \theta', \sigma_1, \sigma_2$ and $1$. Since $\theta +
\theta'= \sigma_1\sigma_2-3$ (equation (\ref{eqn:p>2bis})) we
finally conclude that
\begin{equation}\label{eqn:+}
P^{\mathbb{Z}_3} =\mathbb{Z}[\sigma_1,\sigma_2]\cdot 1 +
\mathbb{Z}[\sigma_1,\sigma_2]\cdot\theta .
\end{equation}
The sign ``$+$'' in formula (\ref{eqn:+}) can be replaced by
``$\oplus$''. Indeed, if $P + Q\theta =0$ for some $P, Q\in
\mathbb{Z}[\sigma_1,\sigma_2]$ then (by interchanging variables
$x$ and $y$) we have $P + Q\theta' =0$ which is possible only if
$P=Q=0$.

By (\ref{eqn:+}) $\theta^2 = P + Q\theta$ for some $P, Q\in
\mathbb{Z}[\sigma_1,\sigma_2]$ which by direct calculation leads
to equation (\ref{eqn:theta-module}) and the isomorphism
(\ref{eqn:theta-ring}). \hfill $\square$

\subsection{The $abc$-lattice $E$ and the $xyz$-lattice $D$}\label{sec:abcxyz}

\begin{figure}[htb]
\centering
\includegraphics[scale=1.00]{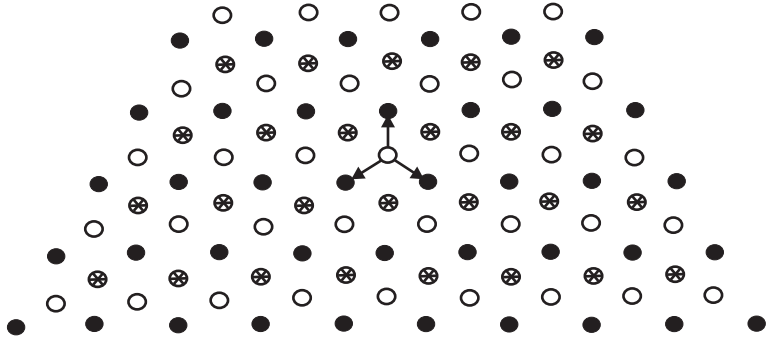}
\caption{The lattice $E$ has three types of dots, black, white,
and asterisks in circles.} \label{fig:3-koseta-1}
\end{figure}

\begin{figure}[htb]
\centering
\includegraphics[scale=0.85]{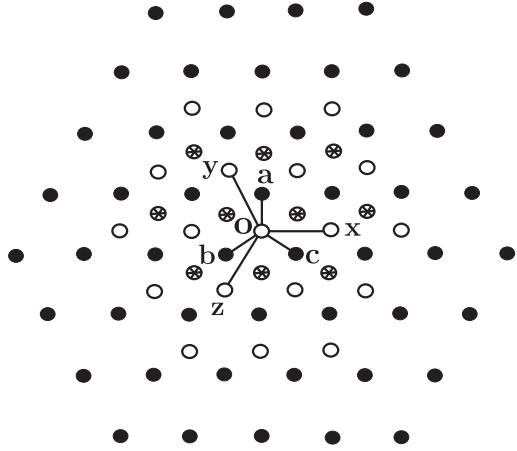}
\caption{Three types of dots correspond to the cosets of the
sublattice $D\subset E$.} \label{fig:tri-koseta-1}
\end{figure}

The `white dot' lattice or the $xyz$-lattice $D$
(Section~\ref{sec:hex-intro}) is generated by vectors
(translations) $t_x,t_y,t_z$ (Figures~\ref{fig:hex-mreza-4} and
\ref{fig:trougaona-mreza-2}). It is convenient to introduce the
$abc$-lattice $E$ as the lattice generated by vectors
$t_a,t_b,t_c$ (Figures~\ref{fig:3-koseta-1} and
\ref{fig:tri-koseta-1}).

The lattice $D$ is a sublattice of $E$ of index $3$. The set of
black dots (Figures~\ref{fig:trougaona-mreza-2},
\ref{fig:3-koseta-1} and \ref{fig:tri-koseta-1}) is clearly one of
the cosets of the quotient lattice $E/D \cong \mathbb{Z}_3$.

\medskip
The fact that $E/D\cong \mathbb{Z}_3$ explains why there are three
types of dots in these images. In order to avoid clutter we will
in subsequent sections continue to draw only black and white dots,
however one should keep in mind the whole of the background
lattice $E$ and the presence of `invisible' dots (circled
asterisks in Figures~\ref{fig:3-koseta-1} and
\ref{fig:tri-koseta-1}).

\noindent Let $\iota : E \rightarrow \mathbb{Z}/ 3\mathbb{Z}$ be
the homomorphism which sends the generators $a,b,c$ to $1\in
\mathbb{Z}$. Then $D = E_0 = \iota^{-1}(0)$ is the set of white
dots and $E_1=\iota^{-1}(1)$ is the set of black dots.

\subsection{The $abc$-ring $Q$ and the $xyz$-ring P}\label{sec:abcxyz-ring}

Let $Q = \mathbb{Z}[a,b,c]/\langle abc-1\rangle$ be the semigroup
ring of the lattice $E$ and let  $Q^{\mathbb{Z}_3} =
(\mathbb{Z}[a,b,c]/\langle abc-1\rangle)^{\mathbb{Z}_3}$ be the
associated ring of $\mathbb{Z}_3$-invariant polynomials. All
structure results that apply to the ring $P =
\mathbb{Z}[x,y,z]/\langle xyz-1\rangle$ apply to the ring $Q$ as
well. In particular there is an isomorphism
(Theorem~\ref{thm:Z3-sim-polynomials }),

\begin{equation}\label{eqn:t-ring}
(\mathbb{Z}[a,b,c]/\langle abc-1 \rangle)^{\mathbb{Z}_3}\cong
\mathbb{Z}[s_1, s_2,t]/\langle \Theta \rangle
\end{equation}
where $\Theta = \Theta(s_1, s_2, t)$ is the polynomial described
in equation (\ref{eqn:theta-module}) and
\begin{equation}\label{eqn:abc-notacija}
s_1=a+b+c \quad s_2=ab+bc+ca \quad t = \Delta(a^2b)=
a^2b+b^2c+c^2a.
\end{equation}
The homomorphism $\iota : E \rightarrow \mathbb{Z}/ 3\mathbb{Z}$
introduced in Section~\ref{sec:abcxyz} allows us to define a
$\mathbb{Z}_3$-grading in the rings $Q$ and $Q^{\mathbb{Z}_3}$.
Indeed, the monomial $m = a^pb^qc^r\in Q$ is graded by its `degree
mod $3$' i.e.\ the mod $3$ class of ${\rm deg}(m)= p+q+r$.

\medskip
The $xyz$-ring $P$, as the group ring of the lattice $D$ of white
dots, corresponds to the elements in the $abc$-ring $Q$ graded by
$0$. Indeed, this follows from the fact
(Figure~\ref{fig:tri-koseta-1}) that $x=ac^2, y=ba^2, z=cb^2$.

Moreover, the $P$-submodule of $Q$ of elements graded by $1$
(monomials with the degree congruent to $1$ mod $3$) is precisely
the submodule generated by the monomials associated with the black
dots.

\medskip
Recall that the $P$-module generated by black dots is precisely
the module $A_{hex}$ from Section~\ref{sec:hex-intro}. This
observation allows us to reduce the ``submodule membership
problem'' in the $P$-module $A_{hex}$ to the corresponding ``ideal
membership problem'' in the ring $Q$. We are primarily interested
in $\mathbb{Z}_3$-invariant polynomials so in the following
section we show how the similar  ``submodule membership problem''
in the $P^{\mathbb{Z}_3}$-module $A_{hex}^{\mathbb{Z}_3}$ can be
reduced to the corresponding ``ideal membership problem'' in the
ring $Q^{\mathbb{Z}_3}$ (Proposition~\ref{prop:membership}).

\subsection{The ring $Q^{\mathbb{Z}_3}$ and $A_{hex}^{\mathbb{Z}_3}$ as a $P^{\mathbb{Z}_3}$-module}
\label{sec:A-hex-Z3-invariant}

The free abelian group $A_{hex}$, generated by all elementary
$2$-cells of the hexagonal lattice $L_{hex}$ (or equivalently all
$0$-dimensional cells of its dual lattice $L_{hex}^\circ$), is a
module over the ring $P = \mathbb{Z}[x,y,z]/\langle xyz-1
\rangle$.

$\Delta$-polynomials already appeared in the description of the
ring $P^{\mathbb{Z}_3}$ of $\mathbb{Z}_3$-invariant polynomials in
$P$ (Figure~\ref{fig:romb-mreza-2}). Following the idea of the
`Newton polygon construction', a polynomial $\Delta(x^py^q) =
 x^py^q+y^pz^q+z^py^q$ is visualized as a triangle with vertices
in the `white dot'-lattice, invariant with respect to the action
of group $\mathbb{Z}_3$.

Similarly the $\Delta$-polynomials $\Delta(a^pb^q) =
 a^pb^q+b^pc^q+c^pa^q$ are, together with $1 = a^0b^0c^0$,
$\mathbb{Z}$-generators of the ring $Q^{\mathbb{Z}_3}$. An
immediate consequence is that $P^{\mathbb{Z}_3}$ is a subring of
$Q^{\mathbb{Z}_3}$.

\begin{figure}[thb]
\centering
\includegraphics[scale=0.80]{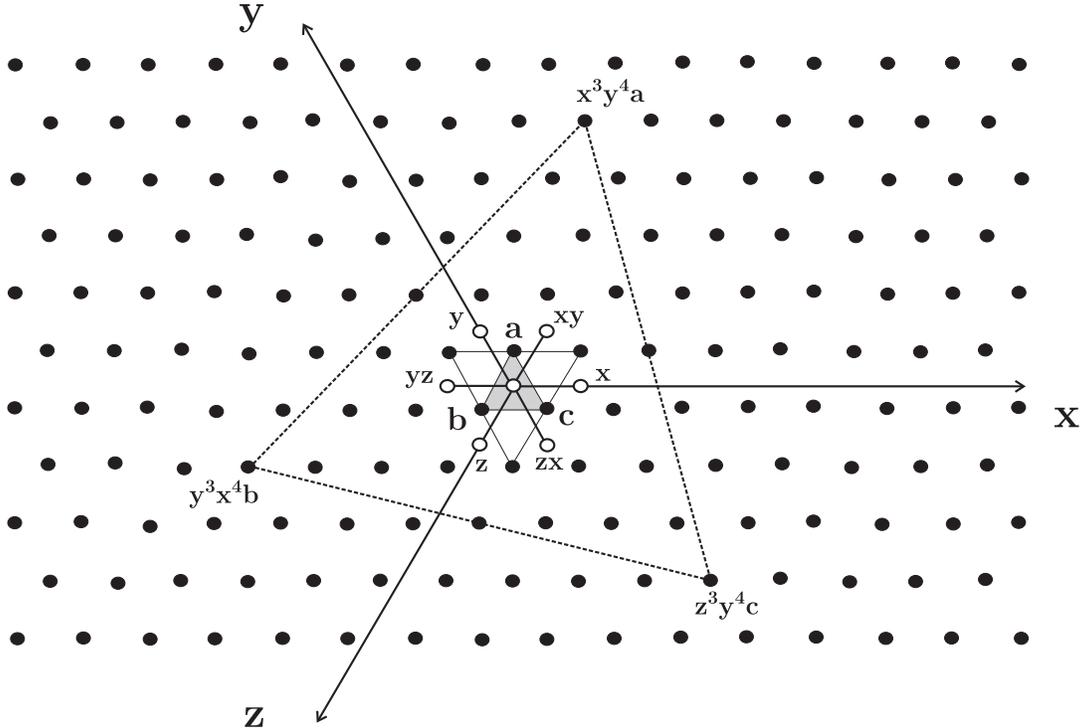}
\caption{$\Delta$-polynomials in the `black dot lattice' generate
the group $A_{hex}^{\mathbb{Z}_3}$.} \label{fig:trougaona-mreza-3}
\end{figure}

Finally, the `black dot' $\Delta$-polynomials $\Delta(x^py^qa)=
x^py^qa+y^pz^qb+z^px^qc$ (Figure~\ref{fig:trougaona-mreza-3}) form
a $\mathbb{Z}$-basis of the group $A_{hex}^{\mathbb{Z}_3}$ of
$\mathbb{Z}_3$-invariant elements of $A_{hex}$.

\medskip
The ring $Q^{\mathbb{Z}_3}$ inherits the $\mathbb{Z}_3$-gradation
from the ring $Q$. Moreover $P^{\mathbb{Z}_3}$ is precisely the
subset of all elements graded by $0\in \mathbb{Z}_3$ while
$A_{hex}^{\mathbb{Z}_3}$ is generated by `black dot triangles'
which are precisely the triangles graded by $1\in \mathbb{Z}_3$.
This characterization is a basis of the following fundamental
proposition.

\begin{prop}\label{prop:membership}
Let $K\subset A_{hex}^{\mathbb{Z}_3}$ be a
$P^{\mathbb{Z}_3}$-submodule of $A_{hex}^{\mathbb{Z}_3}$. Let
$I_K\subset Q^{\mathbb{Z}_3}$ be the ideal in $Q^{\mathbb{Z}_3}$
generated by $K$. Suppose that $p\in A_{hex}^{\mathbb{Z}_3}$.
Then,
\begin{equation}\label{eqn:memb}
p\in K   \quad \iff \quad p\in I_K.
\end{equation}
In other words the `submodule membership problem' is reduced to
the `ideal membership problem' in the ring $Q^{\mathbb{Z}_3}$.
\end{prop}

\medskip\noindent
{\bf Proof:} The implication $p\in K  \Rightarrow p\in I_K$ is
clear. The opposite implication is equally easy. Indeed, if $p =
\alpha_1p_1+\ldots+\alpha_kp_k$ for some elements $p_i\in K$ and
homogeneous (in the sense of the $\mathbb{Z}_3$-gradation)
elements $\alpha_i\in Q^{\mathbb{Z}_3}$ then we can assume that
all $\alpha_i\in P^{\mathbb{Z}_3}$ (the other terms cancel out).
\hfill $\square$

\subsection{Submodule of $A_{hex}^{\mathbb{Z}_3}$ generated by tribones}
\label{sec:submodule}

A three-in-line hexagonal polyomino or a {\em tribone} is a
translate of one of the following three types
(Figure~\ref{fig:trougaona-mreza-3}),

\begin{equation}
\begin{matrix}
T_x & = & x^{-1}+1+x & = & ab^2+1+ac^2\\
T_y & = & y^{-1}+1+y & = & bc^2+1+ba^2\\
T_z & = & z^{-1}+1+z & = & ca^2+1+cb^2
\end{matrix}
\end{equation}
If $A = x^py^qa$ is a `black dot' in the angle $XOY$ then the
three basic tribones centered at the point $A$ are
\begin{equation}\label{eqn:tribone-translates}
T_x(A)= x^py^qaT_x,\, T_y(A)= x^py^qaT_y,\, T_z(A)= x^py^qaT_z.
\end{equation}
For example,
\begin{equation}\label{eqn:tribone-types}
T_x(a) = (x^{-1}+1+x)a \quad T_y(a) = (y^{-1}+1+y)a \quad T_z(a) =
(z^{-1}+1+z)a
\end{equation}
The $\mathbb{Z}_3$-symmetric triplets of tribones, associated to
tribones (\ref{eqn:tribone-translates}), are
\begin{equation}\label{eqn:trib-delta-A}
\Delta(T_x(A)),\quad \Delta(T_y(A)),\quad \Delta(T_z(A))
\end{equation}
where for example,
\begin{equation}\label{eqn:tribone-delta}
\Delta(T_y(A))= \Delta(x^py^qaT_y)=
\Delta(x^py^{q-1}a)+\Delta(x^py^qa)+\Delta(x^py^{q+1}a).
\end{equation}

\begin{theo}\label{thm:tribone-submodule}
The submodule $K_{trib}\subset A_{hex}^{\mathbb{Z}_3}$ of
$\mathbb{Z}_3$-invariant polyominoes (polyhexes) which admit a
signed, symmetric tiling by tribones is generated, as a module
over $P^{\mathbb{Z}_3}$, by the $\mathbb{Z}_3$-symmetric triplets
of tribones,
\begin{equation}\label{eqn:list-delta-tribones}
\Delta(T_x(a)),\, \Delta(T_y(a)),\, \Delta(T_z(a)),\,
\Delta(T_x(ax)),\, \Delta(T_y(ax)),\, \Delta(T_z(ax)).
\end{equation}
\end{theo}

\begin{figure}[!ht]
\centering
\includegraphics[scale=0.80]{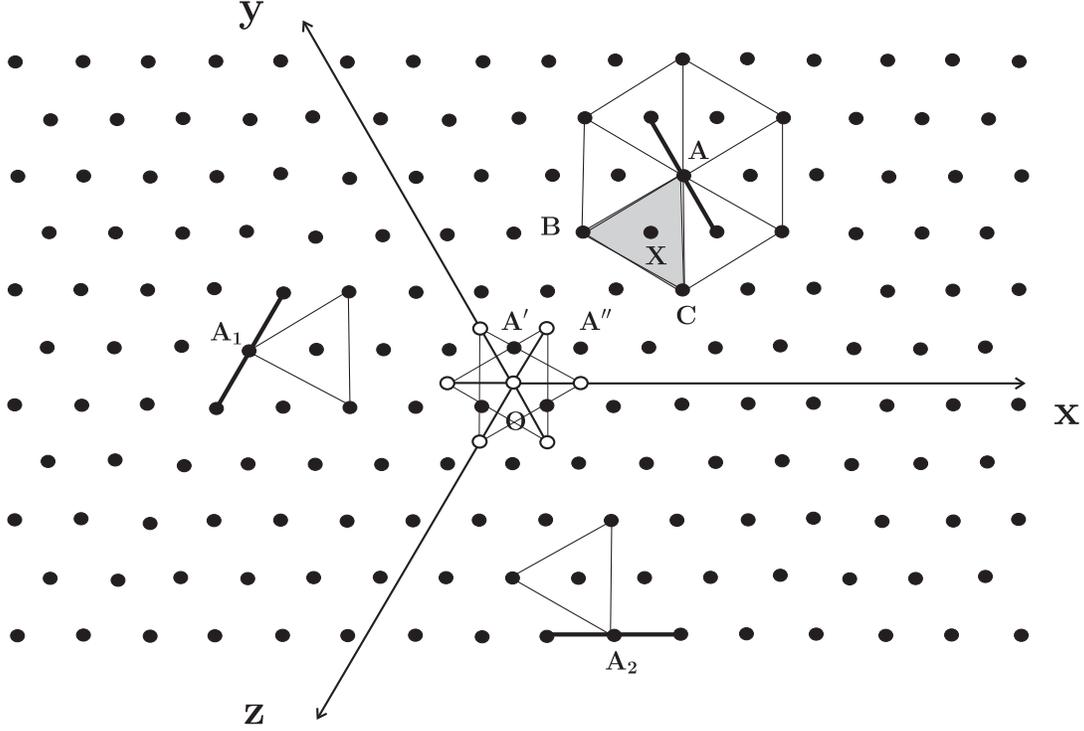}
\caption{$\Delta$-polynomials and the structure of the
$P^\mathbb{Z}_3$-module of $\mathbb{Z}_3$-invariant tiles.}
\label{fig:trougaona-mreza-4}
\end{figure}
\noindent {\bf Proof:} Let $A=Ma = x^py^qa\in A_{hex}$ be a `black
dot' in the angle $\angle XOY$ where $M = x^py^q\in P$ is the
corresponding `white dot' (Figure~\ref{fig:trougaona-mreza-4}).
The three tribones centered at $A$ are $MaT_x, MaT_y, MaT_z$. We
want to show that each of the associated $\mathbb{Z}_3$-symmetric
triplets $\Delta(MaT_x), \Delta(MaT_y), \Delta(MaT_z)$ is in the
$P^{\mathbb{Z}_3}$-module generated by the six elements listed in
the equation (\ref{eqn:list-delta-tribones}).

 \medskip
Figure~\ref{fig:trougaona-mreza-4} depicts the case where
$M=x^4y^3$ and the chosen tribone is $MaT_y=x^4y^3(y^{-1}+1+y)a$.
Using this particular example we describe a general reduction
procedure which allows us to express the triplets $\Delta(MaT_x),
\Delta(MaT_y), \Delta(MaT_z)$ by triplets strictly closer to the
origin $O$.

\medskip
Let $\Delta(x)=x+y+z=\sigma_1$ and $\Delta(xy)=xy+yz+zx=\sigma_2$
be the two basic `white dot' triangles ($\Delta$-polynomials)
depicted at the center of Figure~\ref{fig:trougaona-mreza-4}. By
translating these triangles we surround the point $A$ by six
triangles forming a regular hexagon.

Let us assume that the origin $O$ is {\em not} contained in the
interior of this hexagon. It follows that one of the six triangles
(the shaded triangle $ABC$ in Figure~\ref{fig:trougaona-mreza-4})
has the property that the segment $OA$ intersects the side $BC$
opposite to $A$. (Observe that the origin $O$ can be on the
segment $BC$ and this happens precisely if $A \in \{ya, xya\}$.)

An immediate consequence is that the lengths of all segments $OB,
OC, OX$ are strictly smaller than the length of $OA$. Indeed, by
construction $\angle ACO \geq 60^\circ \geq \angle CAO$, and at
least one of these inequalities is strict. In our case the
triangle $ABC$ is a translate of $\Delta(xy)$ (otherwise we would
use $\Delta(x)$). Since $\Delta(xy)X= A+B+C$ we observe that
\begin{equation}\label{eqn:reduction-1}
\Delta(xy)(y^{-1}+1+y)X = (y^{-1}+1+y)A + (y^{-1}+1+y)B +
(y^{-1}+1+y)C
\end{equation}
Let  $\Delta_Y= \Delta((y^{-1}+1+y)Y)$. By symmetrizing the
equality (\ref{eqn:reduction-1}) and adding (that is by applying
the $\Delta$-operator on both sides of (\ref{eqn:reduction-1})) we
finally obtain,
\begin{equation}\label{eqn:reduction-2}
\Delta(xy)\Delta_X = \Delta_A + \Delta_B + \Delta_C.
\end{equation}
Summarizing, we see that each $\mathbb{Z}_3$-symmetric triplet of
tribones $\Delta_A$ can be expressed in terms of triplets closer
to the origin, provided the hexagon associated to $A$ does not
contain the origin $O$ in its interior.

\medskip
The only remaining possibilities are $A'=a$ and $A''=xa$ which
accounts for the six $\mathbb{Z}_3$-symmetric triplets of tribones
listed in the theorem. \hfill $\square$

\subsection{The ideal $I_{K_{trib}}$}

In this section we express the generating polynomials for the
ideal $I_{K_{trib}}\subset Q^{\mathbb{Z}_3}$ (listed in
(\ref{eqn:list-delta-tribones})) in terms of variables $s_1, s_2,
t$ which appear in the description of the ambient ring
$Q^{\mathbb{Z}_3}$ (equation (\ref{eqn:t-ring}) in
Section~\ref{sec:abcxyz-ring})).

\begin{prop}\label{prop:tribon-ideal}
\begin{equation}
\begin{array}{ccl}
\Delta(T_x(a)) & = & -3 s_1 + 2 s_2^2\\
\Delta(T_y(a)) & = & 3 s_1 - s_2^2 + s_1 t\\
\Delta(T_z(a)) & = & s_1^2 s_2 - s_2^2 - s_1 t\\
\Delta(T_x(ax)) & = & -s_1^2 s_2 + 2 s_2^2 - s_1 t + s_2^2 t\\
\Delta(T_y(ax)) & = & -3 s_1 + s_1^2 s_2 - s_2^2\\
\Delta(T_z(ax)) & = & 3 s_1 - 2 s_1^2 s_2 - s_2^2 + s_1 s_2^3 +
s_1 t - s_2^2 t
\end{array}
\end{equation}
\end{prop}

\medskip\noindent
{\bf Proof:} The proof is by direct calculations which follow the
algorithm described in the proof of
Theorem~\ref{thm:Z3-sim-polynomials }. For example
\[
\Delta(T_x(a))=\Delta(a^2b^2 + a + a^2c^2) =
2\Delta(a^2b^2)+\Delta(a) = 2s_2^2-4s_1+s_1= 2s_2^2-3s_1.
\]
Similarly since
\[ \Delta(a^3b)=\Delta(a)\Delta(a^2b) - \Delta(a^2b^2)-\Delta(a)=
s_1t - (s_2^2-2s_1) -s_1 = s_1t - s_2^2 +s_1
\]
we deduce that,
\[
\Delta(T_y(a))= \Delta(c+a+a^3b) = 2s_1 + \Delta(a^3b) = s_1t -
s_2^2 + 3s_1.
\]

\section{Calculations}
\label{sec:calculations}

\subsection{Auxiliary calculations}\label{sec:auxiliary}

Let $I\subset \mathbb{Z}[x,y]$ be the ideal,
\begin{equation}\label{eqn:aux-ideal}
I = \langle 1+x+x^2, 1+y+y^2, 1+xy+(xy)^2 \rangle .
\end{equation}
If $p-q\in I$ we say that $p$ and $q$ are congruent mod $I$ and
write $p \equiv_I q$. In this section we collect some elementary
congruences mod $I$ which are needed for subsequent calculations.

\begin{lema}\label{lema:mod3}
If $m-n$ is divisible by $3$ then,
\begin{equation}
x^m\equiv_I x^n \qquad y^m\equiv_I y^n \qquad (xy)^m\equiv_I
(xy)^n.
\end{equation}
\end{lema}

\begin{lema}
\begin{equation}\label{eqn:chain}
L_k :=  1+x+\ldots x^{k-1} \equiv_I
 \left\{ \begin{array}{lll}
         0 & \mbox{if\quad $k \equiv 0$ mod $3$}\\
         1 & \mbox{if\quad $k \equiv 1$ mod $3$}\\
         1+x & \mbox{if\quad $k\equiv 2$ mod $3$}.\end{array} \right.
\end{equation}
\end{lema}

\begin{lema}\label{lema:block-square}
\begin{equation}\label{eqn:square}
\square_k := (1+x+\ldots x^{k-1})(1+y+\ldots y^{k-1}) \equiv_I
 \left\{ \begin{array}{lll}
         0 & \mbox{if\quad $k \equiv 0$ mod $3$}\\
         1 & \mbox{if\quad $k \equiv 1$ mod $3$}\\
         x^2y^2 & \mbox{if\quad $k\equiv 2$ mod $3$}.\end{array} \right.
\end{equation}
\end{lema}
\begin{figure}[thb]
\centering
\includegraphics[scale=0.60]{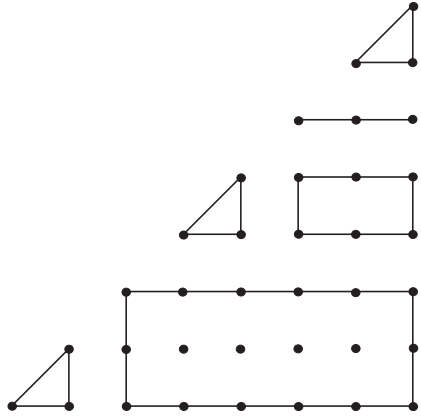}
\caption{Pictorial proof of the relation $\Delta_{3n-1} \equiv_I
n\Delta_2$.} \label{fig:trougao-2}
\end{figure}

\noindent The Newton polygon of the polynomial $\square_k$ is the
square with vertices $\{(0,0), (0, k-1), (k-1, 0), (k-1, k-1)\}$
which has precisely $k$ integer points on each of its sides.
Similarly the polynomial
\begin{equation}
\Delta_k = L_k + xyL_{k-1} + (xy)^2L_{k-2}+\ldots+ (xy)^{k-1}L_1
\end{equation}
collects all monomials associated to the integer points in the
triangle with the vertices $\{(0,0), (k-1,0), (k-1,k-1)\}$
($\Delta_8$ is depicted in Figure~\ref{fig:trougao-2}).

\begin{lema}\label{lema:block-triangle}
\begin{equation}
\begin{array}{lclcc}
\Delta_{3n-1} & \equiv_I & n\Delta_2 &\equiv_I & n(1+x+xy) \\
\Delta_{3n} & \equiv_I &  \Delta_{3n-1} &\equiv_I & n\Delta_2\\
\Delta_{3n+1} & \equiv_I &  n\Delta_2 + 1 &&
\end{array}
\end{equation}
\end{lema}

\subsection{$\Delta$-polynomial of a triangular region}

\begin{figure}[thb]
\centering
\includegraphics[scale=1.20]{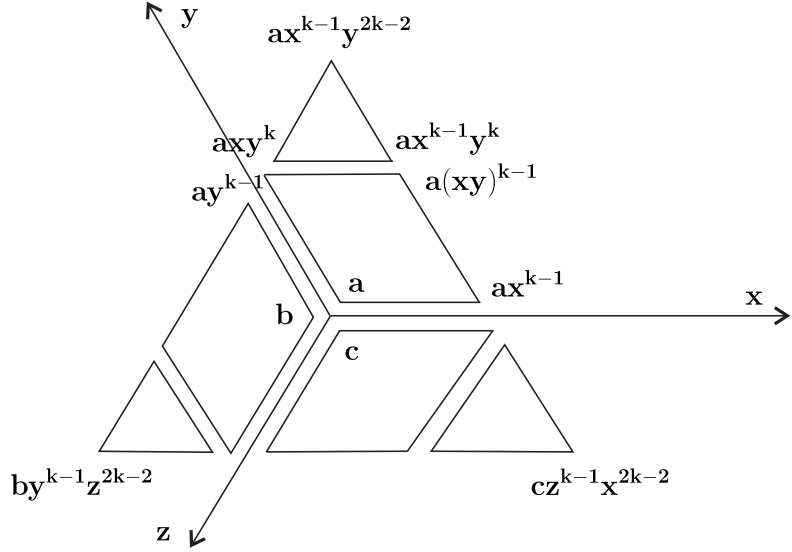}
\caption{Decomposition of $T_N$ into blocks.}
\label{fig:blokovi-3}
\end{figure}

We use the calculations from Section~\ref{sec:auxiliary} to
determine the $\Delta$-polynomial of a triangular region $T_N$
depicted in Figure~\ref{fig:blokovi-3}, where $N=3k-1$ is the
number of black dots on one of the edges of the triangle. By
definition the vertices of the triangular region $T_N$ are black
dots which are associated the monomials,
\[
ax^{k-1}y^{2k-2} \qquad by^{k-1}z^{2k-2} \qquad cz^{k-1}x^{2k-2}.
\]
It is sufficient to determine the $\equiv_I$ class of the
polynomial $A_k$ described as the sum of all monomials in $T_N$
which belong to the cone $xOy$. By inspection we see that $A_k =
B_k + C_k$ where $B_k$ corresponds to the parallelogram with
vertices $a, ax^{k-1}, ay^{k-1}, a(xy)^{k-1}$ and $C_k$ is the
polynomial associated to the triangle with vertices $axy^{k},
ax^{k-1}y^k, ax^{k-1}y^{2k-2}$. Since $B_k = a\cdot\square_k$ and
$C_k = axy^k\Delta_{k-1}$ we can use Lemmas~\ref{lema:mod3},
\ref{lema:block-square} and \ref{lema:block-triangle} to evaluate
$A_k$.

\begin{prop}\label{prop:Ak}
\begin{equation}\label{eqn:Ak}
A_k =
 \left\{ \begin{array}{ll}
         dax\Delta_2 & \mbox{{if}\quad $k=3d$}\\
         a + daxy\Delta_2 & \mbox{if\quad $k=3d+1$}\\
         ax^2y^2+axy^2(d\Delta_2+1) & \mbox{if\quad $k=3d+2$}.\end{array} \right.
\end{equation}
\end{prop}

\medskip\noindent
{\bf Proof:}
\begin{equation}\label{eqn:Bk}
B_k = a\cdot\square_k \equiv_I
 \left\{ \begin{array}{ll}
         0 & \mbox{{if}\quad $k=3d$}\\
         a   & \mbox{if\quad $k=3d+1$}\\
         ax^2y^2  & \mbox{if\quad $k=3d+2$}.\end{array} \right.
\end{equation}
Since,
\begin{equation}\label{eqn:mono}
axy^k \equiv_I
 \left\{ \begin{array}{ll}
         ax & \mbox{{if}\quad $k=3d$}\\
         axy & \mbox{if\quad $k=3d+1$}\\
         axy^2 & \mbox{if\quad $k=3d+2$},\end{array} \right.
\end{equation}
\begin{equation}\label{eqn:Ck}
C_k = axy^k\Delta_{k-1} \equiv_I
 \left\{ \begin{array}{ll}
         dax\Delta_2 & \mbox{{if}\quad $k=3d$}\\
         daxy\Delta_2 & \mbox{if\quad $k=3d+1$}\\
         axy^2(d\Delta_2+1) & \mbox{if\quad $k=3d+2$}.\end{array} \right.
\end{equation}
Since $A_k=B_k + C_k$ the result follows by adding equations
(\ref{eqn:Bk}) and (\ref{eqn:Ck}). \hfill$\square$

\begin{prop}\label{prop:deltaAk}
The $\Delta$-polynomial of $A_k$ is equal to $\Delta(A_k) = P +
dQ$ where,
\begin{flalign*}
\mbox{\rm ($k=3d$)} \hspace{3cm} P=0 &\qquad Q = 3s_1 - 3s_1^2s_2
+ s_1s_2^3 &
\end{flalign*}
\begin{flalign*}
\mbox{\rm ($k=3d+1$)}  \hspace{1cm} P=s_1 &\qquad Q=9s_1 -
6s_1^2s_2 + s_1^3s_2^2 +  4s_1t - 2s_1^2s_2t + s_1t^2 &
\end{flalign*}
\noindent \mbox{\rm ($k=3d+2$)}
\begin{equation*}
\begin{array}{c}
 P= 11 s_1 + s_1^4 - 9 s_1^2 s_2 + 5 s_2^2 + s_1^3 s_2^2 - s_1 s_2^3 +
4 s_1 t -  2 s_1^2 s_2 t + s_2^2 t + s_1 t^2 \\
Q= 24 s_1+s_1^4-11 s_1^2 s_2+s_1^5 s_2-3 s_1^3 s_2^2+4 s_1 s_2^3+8
s_1 t-s_1^4 t-s_1^2 s_2 t+3 s_1 t^2 .
 \end{array}
\end{equation*}
\end{prop}

\medskip\noindent
{\bf Proof:} It follows from Proposition~\ref{prop:Ak} that,
\begin{equation}\label{eqn:Ak-bis}
A_k =
 \left\{ \begin{array}{ll}
         d(a^2c^2+a^3c^4+a^4c^3) & \mbox{{if}\quad $k=3d$}\\
         a+d(a^3c+a^4c^3+a^5c^2) & \mbox{if\quad $k=3d+1$}\\
         (a^5c^2+a^4)+d(a^4+a^5c^2+a^6c) & \mbox{if\quad $k=3d+2$}.\end{array} \right.
\end{equation}
The rest of the proof is by direct calculation, by hand or
preferably by a computer algebra system. \hfill $\square$

\subsection{Gr\"{o}bner basis for the submodule $K_{trib}\subset A_{hex}^{\mathbb{Z}_3}$}

We want to test if the polynomials $\Delta(A_k)$ described in
Proposition~\ref{prop:deltaAk} belong (for different values of
$k$) to the submodule $K_{trib}$ described in
Theorem~\ref{thm:tribone-submodule} (Section~\ref{sec:submodule}).

In light of Proposition~\ref{prop:membership} this question is
reduced to the `ideal membership problem' for the associated ideal
$I_{K_{trib}}$ in the ring $Q^{\mathbb{Z}_3}\cong
\mathbb{Z}[s_1,s_2,t]/\langle \Theta \rangle$ and in turn to the
`ideal membership problem' for the ideal
\begin{equation}
J_{K_{trib}}   :=  I_{K_{trib}} \cup \{\Theta\} \subset
\mathbb{Z}[s_1,s_2,t].
\end{equation}
Here $\Theta$ is the polynomial defined in
Theorem~\ref{thm:Z3-sim-polynomials } (equation
\ref{eqn:theta-ring}) and again, in the context of the $abc$-ring
$Q^{\mathbb{Z}_3}$, in Section~\ref{sec:abcxyz-ring}.

\medskip
In this section we determine the Groebner basis for the ideal
$J_{K_{trib}}$.

\begin{prop}\label{prop:groebner-basis}
The Groebner basis $G=G_{trib}$ of the ideal $J_{K_{trib}} \subset
\mathbb{Z}[s_1,s_2,t]$ with respect to the lexicographic order of
variables $s_1, s_2, t$ is given by the following list of
polynomials:
\begin{equation}\label{eqn:GB-basis}
\begin{array}{cccc}
27+9 t+3 t^2 & -27+t^3  & 9 s_2+3 s_2 t+s_2 t^2 & 3 s_2^2 \\
s_2^2 t & s_2^4  & 3 s_1+s_2^2 & s_2^2+s_1 t \\
s_1 s_2^3 & s_1^2 s_2 & 9+s_1^3+s_2^3+3 t+t^2
\end{array}
\end{equation}
\end{prop}
As the first application of Proposition~\ref{prop:groebner-basis}
we calculate the remainders $\overline{P}^G$ and $\overline{Q}^G$
of polynomials $P$ and $Q$ introduced in
Proposition~\ref{prop:deltaAk}  on division by the Groebner basis
$G=G_{trib}$.

\begin{prop}\label{prop:remainder-P-Q}
Let $P_i$ and $Q_i$ be the polynomials such that to
$\Delta(A_{3d+i}) = P_i + dQ_i$ (Proposition~\ref{prop:deltaAk}).
Then the remainders of these polynomials on division by the
Groebner basis $G=G_{trib}$ are the following:
\begin{equation}\label{eqn:remainders}
\begin{array}{lccc}
\overline{P_0}^G = 0 & & & \overline{Q_0}^G = -s_2^2  \\
\overline{P_1}^G =  s_1 & & & \overline{Q_1}^G = -s_2^2  \\
\overline{P_{2}}^G = -s_1 & & & \overline{Q_{2}}^G = -s_2^2  \\
\end{array}
\end{equation}
\end{prop}

\section{Main results}\label{sec:main}

\begin{theo}\label{thm:main-1}
Let $T_N=T_{3k-1}$ be the $\mathbb{Z}_3$-symmetric triangular
region in the hexagonal lattice depicted in
Figure~\ref{fig:blokovi-3} where $N$ is the number of black dots
(hexagons) on the edge of the triangle. Then $T_N$ admits a
$\mathbb{Z}_3$-symmetric, signed tiling by three-in-line
polyominoes (tribones) if and only if $k = 9r$ for some integer
$r$. The first such triangle is $T_{26}$, in particular the
triangle $T_8$ shown in Figure~\ref{fig:trimino-4} does not have a
$\mathbb{Z}_3$-symmetric, signed tiling by tribones.
\end{theo}

\medskip\noindent
{\bf Proof:} By Proposition~\ref{prop:deltaAk} the polynomial
$\Delta(A_k)$, equal to the sum of all monomials covered by the
triangular region $T_N=T_{3k-1}$, can be expressed (as a
polynomial in variables $s_1,s_2, t$)  as the sum $\Delta(A_k)=
P+dQ$. More explicitly (taking into account the different cases of
Proposition~\ref{prop:deltaAk}) we write $\Delta(A_{3d+i})=
P_i+dQ_i$ where $i\in \{0, 1, 2\}$.

\medskip
By Proposition~\ref{prop:remainder-P-Q} remainders of these
polynomials on division by the Groebner basis $G=G_{trib}$ are
displayed in the table (\ref{eqn:remainders}).

\medskip
The leading terms of the Groebner basis (\ref{eqn:GB-basis}) are
the following,
\begin{equation}\label{eqn:GB-basis-leading-terms}
\begin{array}{ccccccc}
3t^2 & & t^3 & & s_2 t^2 && 3 s_2^2 \\
s_2^2 t && s_2^4  && 3 s_1 && s_1 t \\
s_1 s_2^3 && s_1^2 s_2 && s_1^3
\end{array}
\end{equation}

\medskip By inspection we see that the polynomial $P_i + dQ_i$ can
be reduced to zero if and only $i=0$ and $d=3r$ for some integer
$r$, or in other words if and only if $k=3d=9r$ ($N=27r-1$).
\hfill $\square$

\subsection{The case $N=3k$}\label{sec:case-N3k}

The triangular region $T_N$ in the case $N=3k+1$ has a fixed point
(black dot) with respect the $\mathbb{Z}_3$ action so we focus on
the remaining case $N=3k$.

\medskip\noindent
In this case the convex hull of the set of all black dots in the
intersection of $T_{3k}$ with the angle $xOy$ is the trapeze with
vertices associated the monomials $a, ax^k, ay^{k-1},
ax^{2k-1}y^{k-1}$. We denote the sum of all these monomials by
$A_k$. This trapeze is divided into a rhombus, with vertices at
$a, ax^{k-1}, ay^{k-1}, a(xy)^{k-1}$ and an equilateral triangle
with vertices at $ax^k, ax^ky^{k-1}, ax^{2k-1}y^{k-1}$. The sums
of monomials in these two regions are respectively $B_k$ and
$C_k$, so by definition $A_k=B_k+C_k$.

The mod $I$ class of the polynomial $B_k=a\cdot\square_k$ is as in
the case $N=3k-1$ described by the equality (\ref{eqn:Bk}).
Similarly, $C_k = ax^k \nabla_k$ where $\nabla_k$ is the sum of
all monomials in the triangle with vertices at $1, y^{k-1},
(xy)^{k-1}$. This triangle is obtained from $\Delta_k$ by
interchanging variables $x$ and $y$ so the following lemma is an
immediate consequence of Lemma~\ref{lema:block-triangle}.

\begin{lema}\label{lema:block-nabla}
\begin{equation}
\begin{array}{lclcc}
\nabla_{3n-1} & \equiv_I & n\nabla_2 &\equiv_I & n(1+x+xy) \\
\nabla_{3n} & \equiv_I &  \nabla_{3n-1} &\equiv_I & n\nabla_2\\
\nabla_{3n+1} & \equiv_I &  n\nabla_2 + 1 &&
\end{array}
\end{equation}
\end{lema}
Similarly,
\begin{equation}\label{eqn:mono-2}
ax^k \equiv_I
 \left\{ \begin{array}{ll}
         a & \mbox{{if}\quad $k=3d$}\\
         ax & \mbox{if\quad $k=3d+1$}\\
         ax^2 & \mbox{if\quad $k=3d-1$},\end{array} \right.
\end{equation}
\begin{equation}\label{eqn:Ck-2}
C_k = ax^k\nabla_{k} \equiv_I
 \left\{ \begin{array}{ll}
 da\nabla_2 & \mbox{{if}\quad $k=3d$}\\
 ax(d\nabla_2+1) &\mbox{if\quad $k=3d+1$}\\
 dax^2\nabla_2  &\mbox{if\quad $k=3d-1$}.\end{array} \right.
\end{equation}
and finally,

\begin{prop}\label{prop:Ak-2}
\begin{equation}\label{eqn:Ak-2}
A_k =
 \left\{ \begin{array}{llll}
 da\nabla_2 & = & d(a+a^3b+a^3c) & \mbox{$k=3d$}\\
 a + ax(d\nabla_2+1) & = & (a+a^2c^2)+d(a^2c^2+a^3c+a^4c^3)  & \mbox{$k=3d+1$}\\
 ax^2y^2+dax^2\Delta_2 & = & a^5c^2 + d(a^3c^4+a^4c^3+a^5c^5) & \mbox{$k=3d-1$}.\end{array} \right.
\end{equation}
\end{prop}
The associated $\Delta$-polynomials are given by the following
proposition.

\begin{prop}\label{prop:deltaAk-2}
The $\Delta$-polynomial of $A_k$ is equal to $\Delta(A_k) = P +
dQ$ where,
\begin{flalign*}
\mbox{\rm ($k=3d$)} \hspace{3cm} P=0 &\qquad Q = s_1^2 s_2 - 2
s_2^2 &
\end{flalign*}
\begin{flalign*}
\mbox{\rm ($k=3d+1$)}  \hspace{1cm} P=-s_1 + s_2^2 &\qquad
Q=-s_1^2 s_2 - 2 s_2^2 + s_1 s_2^3 - s_2^2 t &
\end{flalign*}
\noindent \mbox{\rm ($k=3d-1$)}
\begin{equation*}
\begin{array}{c}
 P= 7 s_1 - 5 s_1^2 s_2 + 3 s_2^2 + s_1^3 s_2^2 - s_1 s_2^3 + 4 s_1 t -
 2 s_1^2 s_2 t + s_2^2 t + s_1 t^2 \\
Q= 2 s_1^2 s_2 + 4 s_2^2 - 4 s_1 s_2^3 + s_2^5 .
 \end{array}
\end{equation*}
\end{prop}

\begin{prop}\label{prop:remainder-P-Q-bis}
Let $P_i$ and $Q_i$ be the polynomials such that to
$\Delta(A_{3d+i}) = P_i + dQ_i$
(Proposition~\ref{prop:deltaAk-2}). Then the remainders of these
polynomials on division by the Groebner basis $G=G_{trib}$ are the
following:
\begin{equation}\label{eqn:remainders-2}
\begin{array}{lccc}
\overline{P_0}^G = 0 & & & \overline{Q_0}^G = s_2^2  \\
\overline{P_1}^G =  -s_1 + s_2^2 & & & \overline{Q_1}^G = s_2^2  \\
\overline{P_{-1}}^G = s_1 & & & \overline{Q_{-1}}^G = s_2^2  \\
\end{array}
\end{equation}
\end{prop}

The analysis similar to the proof of Theorem~\ref{thm:main-1}
leads to the following result.
\begin{theo}\label{thm:main-2}
Let $T_N=T_{3k}$ be the $\mathbb{Z}_3$-symmetric triangular region
in the hexagonal lattice which has $3k$ hexagons on one of its
sides. Then $T_N$ admits a $\mathbb{Z}_3$-symmetric, signed tiling
by three-in-line polyominoes (tribones) if and only if $k = 9r$
for some integer $r$. The first such triangle is $T_{27}$.
\end{theo}

\subsection{Examples illustrating Theorems~\ref{thm:main-1} and \ref{thm:main-2}}

Here we describe explicit $\mathbb{Z}_3$-invariant tilings
predicted by Theorems~\ref{thm:main-1} and \ref{thm:main-2}.

\medskip
We begin with the case $N=3k-1= 27r-1$. Figure~\ref{fig:blokovi-3}
describes a decomposition of the triangular region $T_N$ into
triangular and rhombic blocks. The side length of the triangle
with the vertices $axy^{k}, ax^{k-1}y^k, ax^{k-1}y^{2k-2}$ is
$k-1=9r-1$ so by Theorem~\ref{thm:Conv-Laga} this triangle admits
a signed tiling by tribones. The side length of the rhombus with
vertices $a, ax^{k-1}, ay^{k-1}, a(xy)^{k-1}$ is $k=9r$ so it can
be paved by tribones. This signed tiling is extended to the whole
region $T_N$ by rotations through the angle of $120^\circ$ and
$240^\circ$.

\medskip
The case $N=3k= 27r$ is treated similarly. The intersection of
$T_{3k}$ with the angle $xOy$ is the trapeze with vertices at $a,
ax^k, ay^{k-1}, ax^{2k-1}y^{k-1}$ (Section~\ref{sec:case-N3k}).
This trapeze is divided into a rhombus, with vertices at $a,
ax^{k-1}, ay^{k-1}, a(xy)^{k-1}$ and an equilateral triangle with
vertices at $ax^k, ax^ky^{k-1}, ax^{2k-1}y^{k-1}$. Both figures
admit a (signed) tribone tiling. Indeed, the triangle has the side
length $k=9r$ (and the Theorem~\ref{thm:Conv-Laga} applies) while
the side length of the rhombus is as before $k=9r$.

\section{Concluding remarks, examples and questions}\label{sec:concluding}

Our main objective was to illustrate the method of Gr\"{o}bner
bases in the case of equivariant tribone tilings of the hexagonal
lattice. However the method has many other advantages and some of
them were mentioned already in the original paper of Bodini and
Nouvel \cite{BN}.

Perhaps the main reason why this method is so well adapted for
lattice tiling problem is its close connection with the already
developed methods and tools used in lattice geometry.

\subsection{Integer-point transform and Brion's theorem}

Following \cite{Barv} and \cite{BR} we define the {\em
integer-point transform} of a finite subset $K\subset
\mathbb{N}^2$ as the polynomial $f_K = \sum\{x^\alpha y^\beta \mid
(\alpha, \beta)\in \mathbb{N}^2\}$.

Brion's formula, see \cite[Chapter~8]{Barv} and
\cite[Section~9.3]{BR}, is a versatile tool for evaluating the
integer-point transform of convex polytopes. It allows us to
replace a polynomial with a large number of monomials by a very
short expression involving only rational functions. It was tacitly
used throughout the paper for an independent checking of some of
the formulas. More systematic application of these ideas will be
presented in \cite{MMR}.

\subsection{Integer-point enumeration in polyhedra}

Let $P$ be a convex polytope with vertices in $\mathbb{N}^d$ and
let $f_P$ be its integer-point transform. The `discrete volume' of
$Q$, defined as the number of integer points inside $P$, can be
evaluated as the remainder of $f_P$ on division by the ideal
\[
I = \langle x_1-1, x_2-1, \ldots , x_d-1 \rangle.
\]
Let $J\subset \mathbb{Z}[x_1,\ldots, x_d]$ be an ideal, say the
ideal associated to a set $\mathcal{R}$ of prototiles in
$\mathbb{N}^d$. Let $G=G_J$ be the Gr\"{o}bner basis of $J$ with
respect to some term order. It may be tempting to ask (at least
for some carefully chosen ideals $J$) what is the geometric and
combinatorial significance of the remainder $\overline{f}_Q^G$ of
the integer-point transform polynomial $f_Q$ on division by the
Gr\"{o}bner basis $G$.

\bigskip\noindent
{\bf Acknowledgements:} The symbolic algebra computations in the
paper were performed with the aid of {Wolfram Mathematica $9.0$}.


\begin{thebibliography}{9999}


\bibitem{Barv} A. Barvinok. \textit{Integer Points in Polyhedra},
European Mathematical Society, 2008.
\bibitem{BR} M.~Beck and S.~Robins. \textit{Computing the Continuous
Discretely}, Springer 2007.
\bibitem{BW} T. Becker and V. Weispfenning. {\it Gr\"{o}bner Bases}, Springer-Verlag, New York, 1993.
\bibitem{BN} O. Bodini and B. Nouvel. {Z-Tilings of Polyominoes and Standard Basis}, In
\textit{Combinatorial Image Analysis}, Springer (2004), 137--150.
\bibitem{ConLag} J.H. Conway and J.C. Lagarias. Tiling with Polyominoes and
Combinatorial Group Theory, {\it Journal of Combinatorial Theory},
Series A 53 (1990), 183--208.
\bibitem{CLO1}  D. Cox, J. Little and D. O'Shea. \textit{Ideals, Varieties and Algorithms}, Third edition, Springer-Verlag, New York, 2007.
\bibitem{CLO2}  D. Cox, J. Little and D. O'Shea. \textit{Using Algebraic Geometry}, Second edition, Springer-Verlag, New York, 2005.
\bibitem{FT}   D. Fuchs and S. Tabachnikov. \textit{Mathematical Omnibus: Thirty Lectures on Classic
Mathematics},  A.M.S., 2007.
\bibitem{Gol} S.W. Golomb. {\it Polyominoes: Puzzles, Patterns, Problems, and
Packings}, Princeton Univ.\ Press 1996.
\bibitem{MMR} M. Muzika Dizdarevi\'{c}, M. Timotijevi\'{c}, R.T. \v{Z}ivaljevi\'{c}. Signed polyomino tilings, Brion's theorem and Gr\"{o}bner
bases (in preparation).
\bibitem{Reid} M. Reid. {Tile Homotopy Groups}, {\it L'Enseignement Math\`{e}matique} {\textbf 49} (2003), no.1--2,
123--155.
\bibitem{Sturm}  B.~Sturmfels.  \textit{Algorithms in Invariant
Theory}\/ (2nd ed.), Springer 2008.
\bibitem{Thur} W. Thurston. Conway's tiling groups, {\it Amer.\ Math.\
Monthly}\/ 97 (1990), 757--773.

\end{thebibliography}
\end{document}